\input amstex\documentstyle{amsppt}  
\pagewidth{12.5cm}\pageheight{19cm}\magnification\magstep1
\topmatter
\title{Half circles on flag manifolds over a semifield}\endtitle
\author G. Lusztig\endauthor
\address{Department of Mathematics, M.I.T., Cambridge, MA 02139}\endaddress
\thanks{Supported by NSF grant DMS-2153741}\endthanks
\endtopmatter   
\document

\define\mpb{\medpagebreak}

\define\si{\sim}

\define\sqc{\sqcup}

\define\lb{\linebreak}

\define\part{\partial}
\define\emp{\emptyset}

\define\m{\mapsto}
\define\do{\dots}

\define\sub{\subset}    

\define\T{\times}
\define\ti{\tilde}
\define\nl{\newline}
\redefine\i{^{-1}}
\define\fra{\frac}
\define\un{\underline}
\define\ov{\overline}

\define\a{\alpha}
\redefine\b{\beta}
\redefine\c{\chi}
\define\g{\gamma}
\redefine\d{\delta}
\define\e{\epsilon}

\define\io{\iota}

\define\p{\pi}
\define\ph{\phi}

\define\z{\zeta}

\define\ii{\bold i}

\define\kk{\bold k}

\define\HH{\bold H}

\define\NN{\bold N}

\define\RR{\bold R}

\define\ZZ{\bold Z}

\define\cb{\Cal B}

\define\ci{\Cal I}

\define\cu{\Cal U}

\define\cx{\Cal X}

\define\fU{\frak U}

\define\tH{\ti H}

\subhead 1\endsubhead
Let $W$ be a simply laced Weyl group with simple reflections
$\{s_i;i\in I\}$, $I\ne\emp$.
For any semifield $K$ we denote by $\cb_K$
the flag manifold attached to $W,K$ in
\cite{L22, 1.8}. (For $K=\RR_{>0}$, $\cb_K$ was defined in
\cite{L94}.) For any $i\in I$ there
 is a natural free action of $K$ (as a multiplicative group) on
 $\cb_K$  (see no.8). The orbits of this action are called the
 half $i$-circles of $\cb_K$. (This terminology has been used in
\cite{L98} and  \cite{L21, A.2} when $K=\RR_{>0}$, in which case
 $\cb_K$ is an open subset of a real flag manifold and the
 half $i$-circles are the nonempty intersections of $i$-circles
 in the real flag manifold with $\cb_K$.)

We define a graph with $\cb_K$ as set of vertices in which two
points $B\ne B'$ of $\cb_K$ are joined by an edge if there
exists $i\in I$ and a half $i$-circle containing $B$ and $B'$.
Given $k\in\NN$ and $B,B'$ in $\cb_K$ we say that
$B,B'$ are at distance $\le k$ if there exists a sequence
$B=B_1,B_2,\do,B_{k+1}=B'$ in $\cb_K$ such that
$(B_1,B_2),(B_2,B_3),\do,(B_k,B_{k+1})$ are edges of our graph.

 In this paper we shall prove a connectedness property for
the graph $\cb_K$ in the case where

(i) $K$ is a subgroup $\ne\{0\}$ of the additive group $\RR$,
regarded as a semifield with the (new) sum of $a,b$ being
$\min(a,b)$ and the (new) product of $a,b$ being $a+b$.

Let $||:W@>>>\NN$ be the standard length function and let $\nu$
be the maximum value of $||$.

\proclaim{Theorem 2} Assume that $K$ is as in 1(i).

(i) If $B,B'$ are in $\cb_K$ then $B,B'$ are at distance
$\le2\nu-1$.

(ii) In particular, the graph $\cb_K$ is connected.
\endproclaim
This is proved in no.19 as a consequence of a more precise
result (Theorem 13).
In the case $K=\ZZ$, (ii) appears in \cite{L97} with a proof
relying on the theory of canonical bases in quantum groups.

\subhead 3\endsubhead
Here are some examples of semifields other than those in 1(i).

(i) $K$ is a subgroup of the multiplicative group of a field
which is closed under $+$ (such as $\RR_{>0}$ with the usual
$+,\T$).

(ii) $K=\{1\}$ with $1+1=1,1\T1=1$.

If $\kk$ is a totally ordered field, then the subset $\kk_{>0}$
is a semifield as in (i). The analogue of Theorem 2 is still
valid for this
semifield. (See Theorem 20; the proof is similar to that for
$\RR_{>0}$, which appears in \cite{L21, A2}, except that the topological
arguments in {\it loc.cit.} are now replaced by purely algebraic ones.) From this one can
obtain an alternative proof of Theorem 2 (but not of its refinement
in no.13), since the semifield $K$
in no.2 can be viewed as a homomorphic image of a semifield of
form $\kk_{>0}$ (with $\kk$ a field of Puiseux series with real
coefficients with exponents in $K$).

\subhead 4\endsubhead
For $i,j$ in $I$ we set $i.j=-1$ if $i\ne j$,
$s_is_js_i=s_js_is_j$; $i.j=0$ if $i\ne j$, $s_is_j=s_js_i$;
$i.j=2$ if $i=j$.

In the remainder of this paper we fix a semifield $K$.
Following \cite{L19} we consider the monoid $\fU(K)$ with generators
$\{i^a;i\in I,a\in K\}$ and relations

(i) $i^aj^b=j^bi^a$ for $i,j$ in $I$ with $i.j=0$
and $a,b$ in $K$;

(ii) $i^aj^bi^c=j^{a'}i^{b'}j^{c'}$ for $i,j$ in $I$ with
$i.j=-1$ and $a,b,c,a',b',c'$ in $K$ such that
$ab=b'c',a'b'=bc$, $b'=a+c$ (or equivalently $b=a'+c'$);

(iii) $i^ai^b=i^{a+b}$ for $i$ in $I$ and $a,b$ in $K$.

For example, when $K=\{1\}$, the monoid structure of $\fU(\{1\})$
is such that $i^1j^1=j^1i^1$ if $i.j=0$, $i^1j^1i^1=j^1i^1j^1$
if $i.j=-1$ and $i^1i^1=i^1$.

We have a bijection $u\m\un{u}$, $\fU(\{1\})@>\si>>W$, given by
$$\un{i_1^1i_2^1\do i_n^1}=s_{i_1}s_{i_2}\do s_{i_n}$$
for any $i_1,i_2,\do,i_n$ in $I$ such that
$|s_{i_1}s_{i_2}\do s_{i_n}|=n$.
(This bijection is not compatible with the product structures.)

For $u\in\fU(\{1\})$ we set $|u|:=|\un{u}|$. 

A homomorphism of semifields $h:K'@>>>K$ induces a homomorphism
of
monoids $\fU(h):\fU(K')@>>>\fU(K)$ such that $\fU(h)(i^a)=i^{h(a)}$
for any $i\in I$, $a\in K'$.
In particular, the homomorphism of semifields $h_K:K@>>>\{1\}$ given by
$a\m1$ for all $a\in K$ induces a (surjective) homomorphism of monoids
$\fU(h_K):\fU(K)@>>>\fU(\{1\})=W$.
For $w\in W$ let $\fU_w(K)=\fU(h_K)\i(w)$. We have
$\fU(K)=\sqc_{w\in W}\fU_w(K)$.

\subhead 5\endsubhead
For $w\in W$ with $|w|=n$, let $\ci_w$ be the set of sequences
$i_1,i_2,\do,i_n$ in $I$ such that $s_{i_1}s_{i_2}\do s_{i_n}=w$.

For $\ii=(i_1,i_2,\do,i_n)\in\ci_w$ we define $f_\ii(K):K^n@>>>\fU_w(K)$
by
$$f_\ii(K)(a_1,a_2,\do,a_n)=i_1^{a_1}i_2^{a_2}\do i_n^{a_n}.$$

We have the following result.

(a) $f_\ii(K)$ is a bijection.
\nl
The surjectivity of $f_\ii(K)$ is proved as in \cite{L19, 2.9}. We prove
injectivity. When $K$ is as in 3(i) this is proved in
\cite{L19, 2.13}. In the general case we can find a surjective
homomorphism of semifields $h:K'@>>>K$ with $K'$ as in 3(i),
see \cite{BFZ}. We define an equivalence relation on $K'{}^n$ in which the
equivalence classes are the fibres of the map
$h\T h\T\do h:K'{}^n@>>>K^n$. Via the bijection
$f_\ii(K'):K'{}^n@>>>\fU_w(K')$ this becomes an equivalence relation on
$\fU_w(K')$; from the definitions we see that this equivalence relation
is in fact independent of the choice of $\ii$. Let
$\ov{\fU}_w(K')$ be the set of equivalence classes.
Let $\ov{\fU}(K')=\sqc_{w\in W}\ov{\fU}_w(K')$.
From the definitions we see that there is a unique monoid structure
on $\ov{\fU}(K')$ such that the obvious map $\p:\fU(K')@>>>\ov{\fU}(K')$
is a homomorphism of monoids and that the map $i^a\m\p(i^{a'})$ where
$i\in I,a\in K,a'\in h\i(a)\in K'$, extends to a homomorphism of monoids
$\io:\fU(K)@>>>\ov{\fU}(K')$. This restricts to a map
$\io_w:\fU_w(K)@>>>\ov{\fU}_w(K')$. Now $\io_wf_\ii(K)$ coincides with the
bijection $K^n@>>>\ov{\fU}_w(K')$ induced by $f_\ii(K')$. Thus,
$\io_wf_\ii(K)$ is a bijection. It follows that $f_\ii(K)$ is injective.
(In the case where $K=\ZZ$ as in 1(i),
the argument above appears in \cite{L19, 2.13}.)

\subhead 6\endsubhead
Let $w_0$ be the unique element of $W$ such that $|w_0|=\nu$. We
set $\ci=\ci_{w_0}$, $\cu(K)=\fU_{w_0}(K)$. From 5(a) we see
that any $\ii\in\ci$ defines a bijection $f_\ii:K^\nu@>>>\cu(K)$.
This implies that this $\cu(K)$ coincides with the set also
denoted $\cu(K)$ in \cite{L22, 1.2}.
Let $\ph_K:\cu(K)@>>>\cu(K)$ be the involution defined in
\cite{L22, 4.6} in the incarnation \cite{L22} of $\cu(K)$; see
also \cite{L21, 2.1}. (In the case where $K=\RR_{>0}$, $\ph_K$
is implicit in \cite{L94} and is defined in \cite{L97}).
The definition of
$\ph_K$ for general $K$ relies on the properties of $\ph_K$
with $K=\RR_{>0}$ given in \cite{L97}.)

For any $i\in I,a\in K$ there is a well defined map
$T_i(a):\cu(K)@>>>\cu(K)$ such that the following holds:
for any $\ii=(i_1,i_2,\do,i_\nu)\in\ci$ such that $i_1=i$ and any
$(a_1a,a_2,\do,a_\nu)\in K^\nu$ we have
$$T_i(a)(f_\ii(a_1,a_2,\do,a_\nu))=f_\ii(a_1a,a_2,\do,a_\nu).$$
This is proved in \cite{L22, 1.4}. (When $K$ is as in 3(i) this
is proved in \cite{L19, 2.16}; when $K=\ZZ$ as in 1(i) this is
also proved in \cite{L19, 2.16} and is in fact contained in
\cite{L90b} in connection with the theory of canonical bases.)
Note that

$T_i(a)T_i(a')=T_i(aa')$ for $a,a'$ in $K$
\nl
and that $T_i(1)$ is the identity, so that $a\m T_i(a)$ is an
action of $K$ (as a multiplicative group) on $\cu(K)$.

For $i\in I,a\in K$ we have

(a) $T_i(a)\ph_K=\ph_KT_{i^!}(a\i)$
\nl
where $i^!\in I$ is defined by $s_iw_0=w_0s_{i^!}$.
This is proved in
\cite{L22, 1.7} by reduction to the case $K=\RR_{>0}$ which is treated in \cite{L21, 10(a)}.

For any $i\in I$ there is a well defined map $z_i:\cu(K)@>>>K$
such that for any $\ii=(i_1,i_2,\do,i_\nu)\in\ci$ such that
$i_1=i$ and any $(a_1a,a_2,\do,a_\nu)\in K^\nu$ we have
$z_i(f_\ii(a_1,a_2,\do,a_\nu))=a_1$.

From the definition, for any $H\in\cu(K)$ and any $c\in K$ we have

(b) $i^cH=T_i(\fra{c+z_i(H)}{z_i(H)}H$.

\subhead 7\endsubhead
For any function $p:I@>>>K$ there is a well defined monoid isomorphism
$S_p:\fU(K)@>>>\fU(K)$ such that $S_p(i^a)=i^{ap(i)}$ for
$i\in I,a\in K$. This restricts to a bijection
$\cu(K)@>>>\cu(K)$.

For $i\in I,a\in K$ we have

(b) $T_i(a)S_p=S_pT_i(a)$
\nl
as bijections $\cu(K)@>>>\cu(K)$.

We have

(c) $S_p\ph_K=\ph_KS_{p\i}$
\nl
where $p\i(i)=p(i)\i$.
This is proved in
\cite{L22, 1.7} by reduction to the case $K=\RR_{>0}$ which is treated in \cite{L19, 4.3(d)}.

\subhead 8\endsubhead
Let $\cb_K=\{[H,H']\in\cu(K)\T\cu(K);\ph_K(H)=H'\}$,
see \cite{L22, 1.8}. For any $i\in I$ the group $K$ (under
multiplication) acts freely on $\cb_K$ by
$a:[H,H']\m[T_i(a)H,T_{i^!}(a\i)H']$, see 6(a).
The orbits of this action are called half $i$-circles.
Hence $\cb_K$ can be regarded as a graph as in no.1.

We define a graph with $\cu(K)$ as set of vertices in which two
points $H,\tH$ of $\cu(K)$ are joined by an edge if there
exists $i\in I$ and $a\in K$ such that $\tH=T_i(a)H$ (or
equivalently $H=T_i(a\i)\tH$).
Given $k\in\NN$ and $H,\tH$ in $\cb_K$ we say that
$H,\tH$ are at distance $\le k$ if there exists a sequence
$H=H_1,H_2,\do,H_{k+1}=\tH$ in $\cb_K$ such that
$(H_1,H_2),(H_2,H_3),\do,(H_k,H_{k+1})$ are edges of our graph.

We have a bijection $\cu(K)@>\si>>\cb_K$, $H\m[H,\ph_K(H)]$.
This is compatible with the graph structures hence also with
the notion of distance.

\subhead 9\endsubhead
In the remainder of this paper we assume that $K$ is as in
1(a). Let $K_+=K\cap\RR_{\ge0}$. For $w\in W$ with $|w|=n$ let
$\fU_w(K_+)=f_\ii(K)(K_+^n)\sub\fU_w(K)$ where $\ii\in\ci_w$;
from the definitions we see that $\fU_w(K_+)$ is independent of
the choice of $\ii$.
Note that for any $\ii\in\ci_w$, $f_\ii(K):K^n@>>>\fU_w(K)$
restricts to
a bijection $K_+^n@>>>\fU_w(K_+)$. We have the following result.

(a) {\it Let $w\in W$ with $|w|=n$. There is a well defined map
$\c_w:\fU_w(K_+)@>>>\fU(\{1\})$ such that for any
$(i_1,i_2,\do,i_n)\in\ci_w$ and any $(a_1,a_2,\do,a_n)\in K_+^n$
we have
$$\c_w(i_1^{a_1}i_2^{a_2}\do i_n^{a_n})
=i_{j_1}^1i_{j_2}^1\do i_{j_k}^1$$
where $j_1<j_2<\do<j_k$ is the sequence consisting of all
$j\in[1,n]$ such that $a_j=0$.}
\nl
In the case where $K=\ZZ$, this is proved in \cite{L19, \S10}.
The proof in the general case is the same.

We set $\cu(K_+)=\fU_{w_0}(K_+)$.

(The set $\cu(K_+)$ with $K=\ZZ$ appeared in \cite{L90a} as an
indexing set of the canonical basis of the $+$ part of a
quantum group attached to $W$.)

Note that for $\ii=(i_1,i_2,\do,i_\nu)\in\ci$, the element
$i_1^0i_2^0\do i_\nu^0$ of $\cu(K_+)$ is independent of $\ii$; we
denote it by $\HH$. The following result is immediate.

(b) The element $\HH\in\cu(K_+)$ is characterized by the
property that $\c_{w_0}(\HH)=i_\nu^1\do i_2^1i_1^1$.

\subhead 10\endsubhead
Let $i\in I,a\in K$; the bijection $T_i(a):\cu(K)@>>>\cu(K)$
is given by
$$T_i(a)(i_1^{a_1}i_2^{a_2}\do i_\nu^{a_\nu})=
i_1^{a_1+a}i_2^{a_2}\do i_\nu^{a_\nu}$$
for any $(i_1,i_2,\do,i_\nu)\in\ci$ with $i_1=i$ and any
$(a_1,a_2,\do,a_\nu)\in K^\nu$. (Here $a_1+a$ is the usual sum
in $\RR$ which equals the new product of $a_1,a$ in the
semifield $K$.) 
When $a\in K_+$, $T_i(a)$ restricts to an (injective) map
$\cu(K_+)@>>>\cu(K_+)$.

\subhead 11\endsubhead
We fix $(i_1,i_2,\do,i_\nu)\in\ci$ and $k\in\{1,2,\do,\nu\}$.
We have $$(i_k,i_{k+1},\do,i_\nu,i_1^!,i_2^!,i_{k-1}^!)\in\ci.$$

Let $H\in\cu(K_+)$ be such that
$\c_{w_0}(H)=i_{k-1}^1\do i_2^1i_1^1$. We show:

(a) {\it Let $c=z_{i_k}(H)\in K_+$, see no.6. Then
$\tH:=T_{i_k}(-c)H\in\cu(K_+)$ and
$\c_{w_0}(\tH)=i_k^1i_{k-1}^1\do i_2^1i_1^1$. Moreover, we have
$c>0$.}
\nl
Assume that $c=0$. Then $\c_{w_0}(H)=i_k^1u$ for some
$u\in\fU(\{1\})$.
It follows that $i_{k-1}^1\do i_2^1i_1^1=i_k^1u$ hence
$$k=|i_k^1i_{k-1}^1\do i_2^1i_1^1|=|i_k^1i_k^1u|=|i_k^1u|=
|i_{k-1}^1\do i_2^1i_1^1|=k-1,$$
a contradiction. We see that $c>0$. We can find
$c_1,c_2,\do,,c_\nu$ in $K_+$ such that 
$$H=i_k^{c_1}i_{k+1}^{c_2}\do i_\nu^{c_{\nu-k+1}}
(i_1^!)^{c_{\nu-k+2}}(i_2^!)^{c_{\nu-k+3}}\do (i_{k-1}^!)^{c_\nu}.$$
We set
$$H'=i_{k+1}^{c_2}\do i_\nu^{c_{\nu-k+1}}(i_1^!)^{c_{\nu-k+2}}
(i_2^!)^{c_{\nu-k+3}}\do (i_{k-1}^!)^{c_\nu}\in\fU_{s_{i_k}w_0}(K_+).$$
Since $c_1=c>0$, we have $\c_{s_{i_k}w_0}(H')=\c_{w_0}(H)$.
It follows that

(b) $\c_{s_{i_k}w_0}(H')=i_{k-1}^1\do i_2^1i_1^1$.
\nl
We have
$$\tH=i_k^{0}i_{k+1}^{c_2}\do i_\nu^{c_{\nu-k+1}}
(i_1^!)^{c_{\nu-k+2}}(i_2^!)^{c_{\nu-k+3}}\do (i_{k-1}^!)^{c_\nu}.$$
Using (b) we see that 
$\c_{w_0}(\tH)=i_k^1i_{k-1}^1\do i_2^1i_1^1$.
This proves (a).

\subhead 12\endsubhead
We fix $\ii=(i_1,i_2,\do,i_\nu)\in\ci$ and $H_1\in\cu(K_+)$. We
set

(a) $c_1=z_{i_1}(H_1)\in K_+,H_2=T_{i_1}(-c_1)H_1\in\cu(K_+)$,

$c_2=z_{i_2}(H_2)\in K_+,H_3=T_{i_2}(-c_2)H_2\in\cu(K_+)$,

...

$c_\nu=z_{i_\nu}(H_\nu)\in K_+,H_{\nu+1}=T_{i_\nu}(-c_\nu)H_\nu
\in\cu(K_+)$.

We now assume that $c_{w_0}(H_1)=1$. Using 11(a) repeatedly we
see that

$c_1>0,\c_{w_0}(H_2)=i_1^1$,

$c_2>0,\c_{w_0}(H_3)=i_2^1i_1^1$,

...

$c_\nu>0,\c_{w_0}(H_{\nu+1})=i_\nu^1\do i_2^1i_1^1$.

Using 9(b) we see that

(b) $H_{\nu+1}=\HH$.
\nl
It follows that
$$T_{i_\nu}(-c_\nu)\do T_{i_2}(-c_2)T_{i_1}(-c_1)H_1=\HH$$
or equivalently,
$$H_1=T_{i_1}(c_1)T_{i_2}(c_2)\do T_{i_\nu}(c_\nu)\HH$$

\proclaim{Theorem 13}Let $\ii=(i_1,i_2,\do,i_\nu)\in\ci$. Let
$H_1\in\cu(K_+)$. We have
$$H_1=T_{i_1}(c_1)T_{i_2}(c_2)\do T_{i_\nu}(c_\nu)\HH$$
where $c_1,c_2,\do,\c_\nu$ in $K_+$ are defined as in no.12.
\endproclaim
We define $H_2,H_3,\do,H_{\nu+1}$ as in no.12. We show:

(a) $H_{\nu+1}=\HH$.
\nl
Assume first that (a) is known when $K=\RR$. Since $K\sub\RR$ as
a semifield, we can view $\cu(K_+)$ as a submonoid of
$\cu(\RR_{\ge0})$. Let $H_1^\RR$, $H_{\nu+1}^\RR,\HH^\RR$ be the
elements $H_1,H_{\nu+1},\HH$ viewed
as elements of $\cu(\RR_{\ge0})$. From the definitions we see
that
$H_{\nu+1}^\RR$ is obtained from $H_1^\RR$ by the same procedure
(with $K$ replaced by $\RR$) as the one in which $H_{\nu+1}$ is
obtained from $H_1$. Hence using our assumption we have
$H_{\nu+1}^\RR=\HH^\RR$ so that $H_{\nu+1}=\HH$.

Thus it is enough to prove (a) in the case where $K=\RR$. We have
$$H_1=i_1^{a_1}i_2^{a_2}\do i_\nu^{a_\nu},$$
$$H_{\nu+1}=i_1^{b_1}i_2^{b_2}\do i_\nu^{b_\nu},$$
where $(a_1,a_2,\do,a_\nu)\in\RR_{\ge0}^\nu$,
$(b_1,b_2,\do,b_\nu)\in\RR_{\ge0}^\nu$,
are uniquely determined. From the definitions we see that for
any $j$,
$(a_1,a_2,\do,a_\nu)\m b_j$ is a continuous function  
$\RR_{\ge0}^\nu@>>>\RR_{\ge0}$. By 12(b) this function is identically
zero on $\RR_{>0}^\nu$. Being continuous, it is identically zero on
$\RR_{\ge0}^\nu$. This proves (a).
From (a) the theorem follows as in no.12.

\subhead 14\endsubhead
Let $\z)\ii:\cu(K_+)@>>>K_+^\nu$ be the map which to any
$H_1\in\cu(K_+)$ associates $(c_1,c_2,\do,c_\nu)$ as in 1.13.
From no.13 we see that this map is injective. Let $\cx_\ii^K$ be
its
image (a subset of $K_+^\nu$). From the definitions we have

$\cx_\ii^K=\cx_\ii^\RR\cap K_+^\nu$.

\subhead 15\endsubhead
Assume that $I=\{i,j\}$ with $i.j=-1$ so that
$\ii=(i,j,i)\in\ci$.
For $(C,B,A)\in K_+^3$ we have
$$i^Cj^Bi^A=T_i(\g)T_j(\b)T_i(\a)\HH$$
where $\g=C,\b=B+A,\a=B$.
We have $\cx_\ii^K=\{(\g,\b,\a)\in K_+^3;\b\ge\a\}$.

\subhead 16\endsubhead
Assume that $I=\{i,j,k\}$ with $i.j=-1$, $j.k=-1$, $i.k=0$.
We have $\ii=(j,k,i,j,k,i)\in\ci$. For
$(F,E,D,C,B,A)\in K_+^6$ we have
$$j^Fk^Ei^Dj^Ck^Bi^A
=T_j(\ph)T_k(\e)T_i(\d)T_j(\g)T_k(\b)T_i(\a)\HH$$
where
$$\align&\ph=F,\e=B-\min(B,E)+D+C,\d=A-\min(A,D)+E+C,\g=D+E+C,\\&
\b=\min(B,E),\a=\min(A,D).\endalign$$
We have
$$\cx_\ii^K=\{(\ph,\e,\d,\g,\b,\a)\in K_+^6;
\e+\d\ge\g\ge\a+\b,\d\ge\b,\e\ge\a\}.$$

\subhead 17\endsubhead
Assume that $I=\{i,j,k\}$ is as in no.16.
Let $\ii=(i,j,i,k,j,i)$. 

For $(F,E,D,C,B,A)\in K_+^6$ we have
$$i^Fj^Ei^Dk^Cj^Bi^A
=T_i(\ph)T_j(\e)T_i(\d)T_k(\g)T_j(\b)T_i(\a)\HH$$
 where
$$\ph=F,\e=E+D,\d=E,\g=A+B+C,\b=B+C,\a=C.$$
We have
$$\cx_\ii^K=\{(\ph,\e,\d,\g,\b,\a)\in K_+^6;
\e\ge\g,\d\ge\b\ge\a\}.$$

\subhead 18\endsubhead
Let $\ii=(i_1,i_2,\do,i_\nu)\in\ci$. Let $H,\tH$ in $\cu(K)$.
We can find $p:I@>>>K_+$ such that $S_p\tH\in\cu(K_+)$,
$S_pH\in\cu(K_+)$. By Theorem 13 we can find
$(c_1,c_2,\do,c_\nu)\in K_+^\nu$,
$(d_1,d_2,\do,d_\nu)\in K_+^\nu$
such that

(a) $S_pH=T_{i_1}(c_1)T_{i_2}(c_2)\do T_{i_\nu}(c_\nu)\HH$,

(b) $S_p\tH=T_{i_1}(d_1)T_{i_2}(d_2)\do T_{i_\nu}(d_\nu)\HH$.

From (a) we deduce
$$\HH=T_{i_\nu}(-\c_\nu)\do T_{i_2}(-c_2)T_{i_1}(-c_1)S_pH.$$

Introducing this into (b) we obtain
$$S_p\tH=T_{i_1}(d_1)T_{i_2}(d_2)\do T_{i_\nu}(d_\nu)
T_{i_\nu}(-\c_\nu)\do T_{i_2}(-c_2)T_{i_1}(-c_1)S_pH.$$
Using 7(c) we deduce
$$\align&\tH=\\&
T_{i_1}(d_1)T_{i_2}(d_2)\do T_{i_{\nu-1}}(d_{\nu-1})
T_{i_\nu}(d_\nu-\c_\nu)T_{i_{\nu-1}}(-c_{\nu-1})\do
T_{i_2}(-c_2)T_{i_1}(-c_1)H.\tag c\endalign$$

\subhead 19\endsubhead
We prove Theorem 2. It is enough to prove the analogous result
in which the graph $\cb_K$ is replaced by the graph $\cu(K)$
(see no.8).
Let $H\in\cu(K),\tH\in\cu(K)$. Let $(i_1,i_2,\do,i_\nu)$,
$(c_1,c_2,\do,c_\nu)\in K_+^\nu$,
$(d_1,d_2,\do,d_\nu)\in K_+^\nu$ be as in no.18.

We define a sequence

(a) $H_1,H_2,\do,H_{2\nu}$
\nl
in $\cu(K)$ by

$H_1=H$, $H_2=T_{i_1}(-c_1)H_1$, $H_3=T_{i_2}(-c_2)H_2$, ...,
$H_\nu=T_{i_{\nu-1}}(-c_{\nu-1})H_{\nu-1}$,

$H_{\nu+1}=T_{i_\nu}(d_\nu-c_\nu)H_{\nu-1}$,
$H_{\nu+2}=T_{i_{\nu-1}}(d_{\nu-1})H_{\nu+1}$,...,
$H_{2\nu-1}=T_{i_2}(d_2)H_{2\nu-2}$,
$H_{2\nu}=T_{i_1}(d_1)H_{2\nu-1}$,
\nl
From 18(c) we see that (a) is as required in Theorem 2.
This completes the proof.

We have the following result.

\proclaim{Theorem 20} Let $\kk,\kk_{>0}$ be as in no.3. If
$B,B'$ are in $\cb_{\kk_{>0}}$, then $B,B'$ are at distance
$\le2\nu-1$. In particular, the graph $\cb_{\kk_{>0}}$ is
connected.
\endproclaim
It is again enough to prove the analogous result
in which the graph $\cb_{\kk_{>0}}$ is replaced by the graph
$\cu(\kk_{>0})$ (see no.8). The proof (see no.22) uses
the fact that in our case, $\fU(\kk_{>0})$
(and hence $\cu(\kk_{>0})$) is a sub-semigroup of a group
$U(\kk)$ in which $i^a$ is defined for any $a\in\kk$
(not only for $a\in\kk_{>0}$ and we have $i^ai^b=i^{a+b}$ for
any $a,b$ in $\kk$.

\subhead 21\endsubhead
Let $u'\in \cu(\kk_{>0})$, $u''\in\cu(\kk_{>0})$. Let
$i\in I$. We show:

(a) There exists $c\in\kk_{>0}$ such that
$u'i^{-c}\in\cu(\kk_{>0})$, $u''{-c}u''\in\cu(\kk_{>0})$
(equalities in $U(\kk)$).
\nl
We can write
$$u'=i_1^{a_1}i_2^{a_2}\do i_\nu^{a_\nu},
u''=i_1^{b_1}i_2^{b_2}\do i_\nu^{b_\nu}$$
where $i_\nu=i$ and $a_j\in\kk_{>0},b_j\in\kk_{>0}$ for
$j=1,2,\do,\nu$. For $c\in\kk_{>0}$ we have
$$u'i^{-c}=i_1^{a_1}i_2^{a_2}\do i_\nu^{a_\nu-c},
u''i^{-c}=i_1^{b_1}i_2^{b_2}\do i_\nu^{b_\nu-c}.$$
It is enough to show that for some $c\in\kk_{>0}$ we have
$a_\nu-c\in\kk_{>0}$, $b_\nu-c\in\kk_{>0}$. We can assume that
$a_\nu-b_\nu\in\kk_{>0}$. Take $c=(1/2)b_\nu$. We have

$b_\nu-c=c\in\kk_{>0}$, $a_\nu-c=(a_\nu-b_\nu)+(b_\nu-c)
\in\kk_{>0}$.

This proves (a).

\mpb

Let $u'\in\cu(\kk_{>0})$, $u''\in\cu(\kk_{>0}$.
Let $i_1,i_2,\do,i_k$ be a sequence on $I$. We show:

(b) for some $c_1,c_2,\do,c_k$ in $\kk_{>0}$ we have

$u'i_k^{-c_k}i_{k-1}^{-c_{k-1}}\do i_1^{-c_1}\in\cu(\kk_{>0}),
u''i_k^{-c_k}i_{k-1}^{-c_{k-1}}\do i_1^{-c_1}\in\cu(\kk_{>0})$,

(equalities in $U(\kk)$).
\nl
We argue by induction on $k$. If $k=0$ there is nothing to
prove. Assume that $k\ge1$ and the result is known for $k-1$.
Then we can find $c_2,\do,c_k$ in $\kk_{>0}$ such that

$\ti u':=u'i_k^{-c_k}\do i_2^{-c_2}\in\cu(\kk_{>0}),
\ti u'':=u''i_k^{-c_k}\do i_2^{-c_2}\in\cu(\kk_{>0})$.

By (a) we can find $c_1\in\kk_{>0}$ such that

$\ti u'i_1^{-c_1}\in\cu(\kk_{>0}),
\ti u''i_1^{-c_1}\in\cu(\kk_{>0})$,

We see that (b) holds.

We apply (b) with $k=\nu$ and $(i_1,i_2,\do,i_k)\in\ci$. We see
that

(c) for any $u'\in\cu(\kk_{>0})$, $u''\in\cu(\kk_{>0})$ there
exists $u\in\cu(\kk_{>0})$ such that $u'u\i\in\cu(\kk_{>0})$,
$u''u\i\in\cu(\kk_{>0})$. (Here $u'u\i,u''u\i$ are a priori only
in $U(\kk)$.)

\subhead 22\endsubhead
We prove Theorem 20. Let $(i_1,i_2,\do,i_\nu)\in\ci$.
Let $u'\in\cu(\kk_{>0})$, $u''\in\cu(\kk_{>0})$. We choose
$u\in\cu(\kk_{>0})$ as in 21(c). Thus we have

$u'u\i=i_1^{c_1}i_2^{c_2}\do i_\nu^{c_\nu},
u''u\i=i_1^{d_1}i_2^{d_2}\do i_\nu^{d_\nu}$
\nl
(equalities in $U(\kk)$) where
$(c_1,c_2\do,c_\nu)\in\kk_{>0}^\nu$,
$(d_1,d_2\do,d_\nu)\in\kk_{>0}^\nu$.

The following elements are all in $\cu(\kk_{>0})$:

$$\align&u'_\nu=i_1^{c_1}\do i_\nu^{c_\nu}u,
\do,u'_2=i_{\nu-1}^{c_{\nu-1}}i_\nu^{c_\nu}u,
u'_1=i_\nu^{c_\nu}u,u'_0=u''_0=u,\\&
u''_1=i_\nu^{d_\nu}u,
u''_2=i_{\nu-1}^{d_{\nu-1}}i_\nu^{d_\nu}u,\do,
u''_\nu=i_1^{d_1}\do i_\nu^{d_\nu}u.\endalign$$

Any two consecutive terms of this sequence form an edge of our
graph (we use 6(b)); moreover, even after we drop
$u'_0=u''_0=u$ from the sequence, any two consecutive terms of
the resulting sequence form an edge of our graph: we use that
$u'_1=i_\nu^{c_\nu d_\nu\i}u''_1$.
This completes the proof.

\subhead 23\endsubhead
Let $\cb$ be the set of Borel subgroups defined over $\kk$
of a split, connected reductive group over $\kk$ with Weyl group
$W$. We can identify $\cb_{\kk_{>0}}$ with a subset of $\cb$ as
in \cite{L22, 4.8}. We define a graph structure on
$\cb_{\kk_{>0}}$ in which two Borel subgroups $B,B'$ in
$\cb_{\kk_{>0}}$ are joined if their relative position (a Weyl
group element) is $s_i$ for some $i\in I$.
By arguments similar to those in \cite{L21, A2} we see that this
graph structure is the same as that in no.1. It follows that
this graph is connected.

\enddocument